\documentclass[12pt,a4paper]{article}
\usepackage{amsfonts,amsmath,amsthm,amssymb,stmaryrd}
\newtheorem{theorem}{Theorem}[section]
\newtheorem{lemma}[theorem]{Lemma}
\newtheorem{Proposition}[theorem]{Proposition}

\numberwithin{equation}{section}

\usepackage[left=1 in, right=1 in,top=1 in, bottom=1 in]{geometry}

\begin{document}
  \title{{Critical Magnetic Number in the  MHD\\ Rayleigh-Taylor instability}\thanks{Partially supported by National Natural Science Foundation of China-NSAF (No. 10976026).}}
  \author{Yanjin Wang\\{\it\small School of Mathematical Sciences, Xiamen University,
 Fujian 361005, China}\\{\it\small
 Email: yanjin$\_$wang@xmu.edu.cn}}
  \maketitle
    \begin{abstract}{\small
We reformulate in Lagrangian coordinates the two-phase free boundary problem for the equations of Magnetohydrodynamics in a infinite slab, which is incompressible, viscous and of zero resistivity, as one for the Navier-Stokes equations with a force term  induced by the fluid flow map. We study the stabilized effect of the magnetic field for the linearized equations around the steady-state solution by assuming that the upper fluid is heavier than the lower fluid, $i. e.$, the linear Rayleigh-Taylor instability.  We identity the
 critical magnetic number $|B|_c$ by a variational problem. For the cases $(i)$ the magnetic number $\bar{B}$ is vertical in 2D or 3D; $(ii)$ $\bar{B}$ is horizontal in 2D, we prove that the linear system is stable when  $|\bar{B}|\ge |B|_c$ and is unstable when $|\bar{B}|<|B|_c$. Moreover, for $|\bar{B}|<|B|_c$ the vertical $\bar{B}$ stabilizes the low frequency interval while the horizontal $\bar{B}$ stabilizes the high frequency interval, and the growth rate of growing modes is bounded.}\smallskip
\end{abstract}

\section{Formulation}
\subsection{Formulation in Eulerian coordinates}
We consider the two-phase free boundary problem for the equations of Magnetohydrodynamics (MHD)
 within the infinite slab
$\Omega=\mathbb{R}^{d-1}\times(-1,1)\subset \mathbb{R}^{d},\ d=2\hbox{ or }3$ is the dimension, and
for time $t\ge0$. The fluids are separated by a moving free
interface $\Sigma(t)$ that extends to infinity in every horizontal
direction. The interface divides $\Omega$ into two time-dependent,
disjoint, open subsets $\Omega_\pm(t)$ so that $\Omega=
\Omega_+(t)\cup\Omega_-(t)\cup\Sigma(t)$ and
$\Sigma(t)=\bar{\Omega}_+(t)\cap\bar{\Omega}_-(t)$.  The
motions of the fluids are driven by the constant gravitational field
along $e_d-$the $x_d$ direction, $G=(0,0,-g)$ with $g>0$ and the Lorentz force induced by the magnetic fields. The two
 fluids are described by their  velocity, pressure and magnetic field
 functions, which are given for each $t\ge0$ by, respectively,
\begin{equation}(u_\pm, \bar{p}_\pm, h_\pm)(t,\cdot):\Omega_\pm(t)\rightarrow
(\mathbb{R}^+,\mathbb{R}^d,\mathbb{R},\mathbb{R}^d).\end{equation}
We shall assume that at a given time $t\ge0$ these functions have well-defined traces onto $\Sigma(t)$.

The fluids under consideration are incompressible, viscous and of zero resistivity, hence
for $t>0$ and $(x',x_d)\in \Omega_\pm(t)$ the fluids satisfy the following magnetohydrodynamic equations:

$$\left\{\begin{array}{ll}\partial_t(\rho_\pm u_\pm) + {\rm div}(\rho_\pm u_\pm\otimes u_\pm)+{\rm
div}S_\pm=-g\rho_\pm e_d,
\\ {\rm div}u_\pm=0,\\\partial_t h_\pm +{\rm div}(u_\pm\otimes h_\pm)-{\rm div}(h_\pm\otimes u_\pm)=0,
\\ {\rm div}h_\pm=0, \end{array}\right.\eqno(1.2)$$where we define the stress tensor consisting of fluid part and magnetic part
by\setcounter{equation}{2}
\begin{equation}
S_{\pm}=-\mu_{\pm}(\nabla u_\pm+\nabla u_\pm^T)+ \bar{p}_\pm I +\frac{|h_\pm|^2}{2}I  -h_\pm\otimes h_\pm.
\end{equation}
Hereafter the superscribe $T$ means the transposition and $I$ is the
$d\times d$ identity matrix. The positive constants $\rho_\pm$ and $\mu_\pm$ denote the densities and viscosities of the respective fluids.

Now we prescribe the boundary conditions at the fixed boundaries and the jump conditions at the free interface. First for the the fluid equations $(1.2)_1$, due to the viscosity we assume that the velocity is continuous across the free interface and enforce the no-slip condition at the fixed upper and lower boundaries. Since we do not take into account the surface tension, it is standard to assume that the normal
stress is continuous across
the free interface (cf. \cite{C,WL}). Therefore, we impose the boundary conditions at the fixed boundaries
\begin{eqnarray}&&
 u_+(t,x',-1)=u_-(t,x',1)=0,\hbox{ for all }t\ge 0,\ x'\in
 \mathbb{R}^{d-1},\end{eqnarray}
and impose
the jump conditions at the free interface
\begin{eqnarray}&&
[u]\mid_{\Sigma(t)}=0,
\\&& [S\nu]\mid_{\Sigma(t)}=0,\end{eqnarray}
where we have written the normal vector to $\Sigma(t)$  as
$\nu$,  $f|_{\Sigma(t)}$ for the trace of a quantity $f$ on $\Sigma(t)$ and denote the interfacial jump by
\begin{equation}[ f]\mid_{\Sigma(t)}:=f_+|_{\Sigma(t)}-f_-|_{\Sigma(t)}.
\end{equation}For the magnetic equations $(1.2)_3$, since the fluids are of zero resistivity this is a free transport equation along the flow, and hence the Dirichlet boundary condition on the velocity at the fixed boundary prevents the necessity of prescribing boundary condition on the magnetic fields. On the other hand,
 due to the divergence-free $(1.2)_4$ and also physically,
 we shall assume that the normal component of  magnetic field is continuous across the free interface (cf. \cite{C,S})
 \begin{eqnarray}
[h\cdot \nu]|_{\Sigma(t)}=0.\end{eqnarray}
However, we will show explicitly   that the divergence-free of $h_\pm$ and the jump condition (1.8) can be verified if they hold initially. Therefore, the conditions $(1.2)_4$,  (1.8) are transformed to  the compatibility conditions assumed on the initial magnetic field.

The motion of the free interface is coupled to the evolution
equations for the fluids (1.2) by requiring that the boundary be advected
with the fluids. More precisely, if $V(t,x)\in \mathbb{R}^d$ denotes
the normal velocity of the boundary at $x\in\Sigma(t)$, then
\begin{equation}V(t,x) = (u(t,x)\cdot\nu(t,x))\nu(t,x).\end{equation}
Here $u(t,x)$ is the common trace of $u_\pm(t,x)$ onto  $\Sigma(t)$
and these traces agree because of the  jump condition (1.5),
which also implies that there is no possibility of the fluids
slipping past each other along $\Sigma(t)$.

To complete the statement of the problem, we must specify initial
conditions. We give the initial interface $\Sigma(0) = \Sigma_0$,
which yields the open sets $\Omega_\pm(0)$ on which we specify the
initial data for the velocity and magnetic
field
\begin{equation}(u_\pm,h_\pm)(0,\cdot):\Omega_\pm(0)\rightarrow
(\mathbb{R}^d,\mathbb{R}^d).\end{equation}

To simply the equations we introduce the indicator function $\chi$ and denote
$$\left\{\begin{array}{ll}u=u_+{\chi_{\Omega_+}}+u_-{\chi_{\Omega_-}},\quad h=h_+{\chi_{\Omega_+}}+h_-{\chi_{\Omega_-}}, \quad \bar{p}=\bar{p}_+{\chi_{\Omega_+}}+\bar{p}_-{\chi_{\Omega_-}},\\\rho=\rho_+{\chi_{\Omega_+}}+\rho_-{\chi_{\Omega_-}},\quad \mu=\mu_+{\chi_{\Omega_+}}+\mu_-{\chi_{\Omega_-}}
,\end{array}\right.\eqno(1.11)$$\setcounter{equation}{11}
and also define the modified pressure by
\begin{eqnarray}
p=\bar{p}+\frac{|h|^2}{2}+g\rho x_d.\end{eqnarray}
Hence the equations (1.2) are replaced by
$$\left\{\begin{array}{ll}\rho\partial_t u+ \rho u\cdot\nabla u-\mu\Delta u+\nabla p=h\cdot\nabla h,\\
\partial_t h +u\cdot\nabla h-h\cdot\nabla u=0,
\\{\rm div}u={\rm div}h=0, \end{array}\right.\eqno(1.13)$$\setcounter{equation}{13}
and the jump condition (1.6) becomes, setting $[\rho]=\rho_+-\rho_-$,
\begin{equation}\left[\left(p I-\mu(\nabla u+\nabla u^T)\right)\nu\right]\mid_{\Sigma(t)}=g[\rho]x_d\nu+h\cdot\nu\left[h \right]\mid_{\Sigma(t)}.\end{equation}

\subsection{Formulation in Lagrangian coordinates}

Since the movement of the free interface $\Sigma(t)$ and the
subsequent change of the domains $\Omega_\pm(t)$ in Eulerian
coordinates result severe mathematical difficulties,  we switch our
analysis to Lagrangian coordinates so that the interface and the
domains stay fixed in time. To this end we define the fixed
Lagrangian domains $\Omega_+=\mathbb{R}^{d-1}\times(0,1)$ and
$\Omega_-=\mathbb{R}^{d-1}\times(-1,0)$. We assume that there exist
invertible mappings
\begin{equation}\eta_\pm^0 :\Omega_\pm \rightarrow
\Omega_\pm(0), \end{equation} which are continuous across $\{x_d=0\}$
so that $\Sigma_0=\eta^0_\pm(\{x_d=0\}),\
\{x_d=1\}=\eta^0_+(\{x_d=1\})$, and $\{x_d=-1\}=\eta^0_-(\{x_d=-1\})
$. The first condition means that $\Sigma_0$ is parameterized by the
either of the mappings $\eta_\pm$ restricted to $\{x_2 = 0\}$ (which
one is irrelevant since they are continuous across the interface),
and the latter two conditions mean that $\eta_\pm$ map the fixed
upper and lower boundaries into themselves. Define the flow maps
$\eta_\pm$ as the solution to
$$\left\{\begin{array}{ll}
\partial_t\eta_\pm(t,x)=u_\pm(t,\eta_\pm(t,x)),\\
\eta_\pm(0,x)=\eta_\pm^0(x).
\end{array}\right.\eqno(1.16)$$
We think of the Eulerian coordinates as $(t,y)$ with $y =
\eta(t,x)$, whereas we think of Lagrangian coordinates as the fixed
$(t,x)\in \mathbb{R}^+\times\Omega$. In order to switch back and
forth from Lagrangian to Eulerian coordinates we assume that that
$\eta_\pm(t,\cdot)$ are invertible and  $\Omega_\pm(t) = \eta_\pm(
t,\Omega_\pm)$, and since $u_\pm$ and $\eta_\pm^0$ are all
continuous across $\{x_d = 0\}$, we have $\Sigma(t) =
\eta_\pm(t,\{x_d = 0\})$. In other words, the Eulerian domains of
upper and lower fluids are the image of $\Omega_\pm$ under the
mappings $\eta_\pm$ and the free interface is the image of $\{x_d =
0\}$ under the mapping $\eta_\pm(t,\cdot)$.

 Setting $\eta=\chi_+\eta_++\chi_-\eta_-$, we define the Lagrangian unknowns\setcounter{equation}{16}
\begin{equation}(v,b,q)(t,x)=(u,h,p)(t,\eta(t,x)),\quad (t,x)\in \mathbb{R}^+\times\Omega.\end{equation}
 Defining the matrix $A$ via $A^T = (D\eta)^{-1}$, then
in Lagrangian coordinates the evolution equations for $\eta,v,b,q$
are, writing $\partial_j =
\partial/\partial_{x_j},$\setcounter{equation}{17}
$$\left\{\begin{array}{ll}
\partial_t\eta_i=v_i,
 \\\rho\partial_tv_i+A_{jk}\partial_kT_{ij}=b_jA_{jk}\partial_kb^i,
 \\ A_{jk}\partial_kv_j =0,
 \\ \partial_tb_i=b_jA_{jk}\partial_kv_i,
\\  A_{jk}\partial_kb_j=0,\end{array}\right.\eqno(1.18)$$
where the stress tensor of fluid part in Lagrangian coordinates,  $T(v,q)$, is given
by\setcounter{equation}{18}
\begin{equation}
T_{ij}=qI_{ij}-\mu(A_{jk}\partial_kv_i+A_{ik}\partial_kv_j).
\end{equation}
Here we have written $I_{ij}$ for $i, j$ component of the identity
matrix $I$ and we have employed the Einstein convention of summing
over repeated indices.

To write the jump conditions, for a quantity $f = f_\pm$, we define the interfacial jump as
\begin{equation}\llbracket f\rrbracket:=f_+|_{\{x_d=0\}}-f_-|_{\{x_d=0\}}.
\end{equation}
Then the jump conditions in Lagrangian coordinates are
\begin{equation}\llbracket v\rrbracket=0, \quad \llbracket b_jn_j\rrbracket=0,\quad
\llbracket T_{ij}n_j\rrbracket=g[\rho]\eta_dn^i+ \llbracket b_i\rrbracket b_jn_j,
\end{equation}
where the unit normal to the interface $\Sigma(t)$, $i.e.,$ $n=\nu(\eta)$ can be represented by
\begin{equation}n=\left.\frac{Ae_d}{|Ae_d|}\right|_{\{x_d=0\}}.
\end{equation}
Finally, we require the no-slip boundary condition
\begin{eqnarray}&&
 v_+(t,x',-1)=v_-(t,x',1)=0,\hbox{ for all }t\ge 0,\ x'\in
 \mathbb{R}^{d-1}.\end{eqnarray}

\subsection{Reformulation}
One purpose of this paper is to reformulate the free boundary
problem (1.18), (1.21), (1.23) as the Navier-Stokes equations with a
force term  induced by the fluid flow map. We want to eliminate $b$
by expressing it in terms of $\eta$. Indeed, applying $A_{il}$ to
$(1.18)_4$, we have
\begin{equation}A_{il}\partial_tb_i=b_jA_{jk}\partial_kv_iA_{il}=b_jA_{jk}\partial_t(\partial_k\eta_i)A_{il}
=-b_jA_{jk}\partial_k\eta_i \partial_tA_{il}=-b_i\partial_tA_{il}.\nonumber  \end{equation}
This implies that $\partial_t (A_{jl}b_j)=0$ and hence,
\begin{eqnarray}
&&A_{jl}b_j=A_{jl}^0b_j^0,\\&&
b_i =\partial_l\eta^iA^0_{jl}b^0_j.\end{eqnarray}
Hereafter, the superscript $0$ means the initial value.

By the expression (1.25), we first check the divergence of $b,\ i.e.,$ $(1.18)_5$. We make use of the geometric identities
\begin{equation}J=J^0\hbox{ and }\partial_k(JA_{ik})=0.\end{equation}
Hence, applying $A_{ik}\partial_k$ to (1.25) we have
\begin{equation}
A_{ik}\partial_kb_i=\frac{J}{J^0}A_{ik}\partial_k(\partial_l\eta^iA^0_{jl}b^0_j)
=\frac{1}{J^0}\partial_k(JA_{ik}\partial_l\eta^iA^0_{jl}b^0_j)
=\frac{1}{J^0}\partial_k(J^0A^0_{jk}b^0_j)=A^0_{jk}\partial_k b^0_j.\end{equation}
Next we check the jump  $\llbracket b_jn_j\rrbracket$. We  recall the expression (1.22) of unit normal to the interface $\Sigma(t)$. It is easy to verify that $Ae_d$ is continuous across the free interface. Hence we have
\begin{equation}\llbracket b_jn_j\rrbracket=\llbracket\partial_l\eta^jA^0_{kl}b^0_k\cdot{A_{jd}}\rrbracket\frac{1}{|Ae_d|}=\llbracket A^0_{kd}b^0_k\rrbracket\frac{1}{|Ae_d|}=\llbracket b_j^0n_j^0\rrbracket \frac{|A^0e_d|}{|Ae_d|}.
\end{equation}
Hence, if we assume the compatibility conditions on the initial data
\begin{equation}A^0_{jk}\partial_k b^0_j=0,\quad \llbracket b_j^0n_j^0\rrbracket=0,\end{equation}
then from (1.27), (1.28), we have
\begin{equation}A_{jk}\partial_k b_j=0,\quad \llbracket b_jn_j\rrbracket=0.\end{equation}

Moreover, for simplify of notations, we assume that
\begin{equation}A^0_{ml}b^0_m=\bar{B}_l\hbox{ with }\bar{B} \hbox{ is a constant vector}.\end{equation}
We remark that the class of the pairs of the data $\eta^0, b^0$ that satisfy the constraints (1.29), (1.31) is quite large enough. For example, we chose $\eta^0=Id$ and $b^0=const$, then by (1.24), (1.27), (1.28), any pair of data $\eta, b$ which is transported by the flow will satisfy (1.29), (1.31).

Now we represent the Lorentz force term by, since (1.24), (1.25), (1.30), (1.31),
\begin{equation}b_jA_{jk}\partial_kb^i=\partial_l\eta^jA^0_{ml}b^0_mA_{jk}\partial_k(\partial_r\eta^iA^0_{sr}b^0_s)
=A^0_{mk}b^0_m\partial_k(\partial_r\eta^iA^0_{sr}b^0_s)=\bar{B}_l\bar{B}_m\partial_{lm}^2\eta^i.\end{equation}
Hence the equations (1.18) becomes  a Navier-Stokes system with the force term induced by the flow map $\eta$:
$$\left\{\begin{array}{ll}
\partial_t\eta_i=v_i,
 \\\rho\partial_tv_i+A_{jk}\partial_kT_{ij}-\bar{B}_l\bar{B}_m\partial_{lm}^2\eta^i=0
,\\ A_{jk}\partial_kv_j=0,\end{array}\right.\eqno(1.33)$$where the magnetic number $\bar{B}$ can be regarded as a vector parameter. Accordingly, the jump conditions (1.21) become\setcounter{equation}{33}
\begin{equation}\llbracket v\rrbracket=0, \quad
\llbracket T_{ij}n_j\rrbracket=g[\rho]\eta_dn^i+\bar{B}_l \bar{B}_m \llbracket\partial_l\eta^i\rrbracket\partial_m\eta^jn^j,
\end{equation}
Note that  we implicitly admit that $\bar{B}_m\partial_m\eta^jn^j$ is continuous across $\{x_d=0\}$, this follows from the assumptions (1.29), (1.31). Finally, we require the  boundary condition (1.23).

\subsection{Linearization around the steady state}
 The system (1.33), (1.34), (1.23) admits the steady solution with $v=0,\ \eta=Id, ,\ q=const$ with the interface given by
 $\eta(\{x_d=0\})=\{x_d=0\}$ and hence $n=e_d,\ A=I$. Now we linearize the equations (1.33) around such a steady-state
 solution, the resulting linearized equations are
\setcounter{equation}{34}
$$\left\{\begin{array}{ll}\partial_t\eta=v,\\
\rho\partial_tv +\nabla q-\mu\Delta v-\bar{B}_l\bar{B}_m\partial_{lm}^2\eta=0,
\\{\rm div}v=0.\end{array}\right.\eqno(1.35)$$
The corresponding linearized jump
conditions are\setcounter{equation}{35}
\begin{eqnarray}&&\llbracket v\rrbracket=0,\quad\llbracket -\mu (Dv+Dv^T) + qI\rrbracket e_d =g[\rho]\eta_de_d+\bar{B}_d\bar{B}_l\llbracket\partial_{l}\eta\rrbracket,
\end{eqnarray}
while the boundary conditions are
\begin{equation}v_-(t,x',-1)=v_+(t,x',1)=0.
\end{equation}

The main purpose of this paper is to study the stabilized effect of
magnetic field on the Rayleigh-Taylor problem, hence we assume that
the upper fluid is heavier than the lower fluid, $i.e.,$
\begin{equation}\rho_+>\rho_-\Longleftrightarrow[\rho]>0.
\end{equation}

\subsection{Normal mode ansatz}
In the fluid stability analysis it is standard to study the normal mode solutions of the linearized system (cf. \cite{C}).
 To begin, we assume an ansatz
\begin{equation}v(t,x)=w(x){\rm e}^{\lambda t},\ q(t,x)=\tilde{q}(x){\rm e}^{\lambda t},
\ \eta(t,x)=\tilde{\eta}(x){\rm e}^{\lambda t},\end{equation}for some $\lambda>0$. Substituting
this ansatz into (1.35), eliminating $\tilde{\eta}$ by using the first equation, we arrive at
the time-invariant system for $w=(w_1,\dots,w_d)$ and $\tilde{q}$:
$$\left\{\begin{array}{ll} \lambda\rho w +\nabla \tilde{q} -\mu\Delta w-\lambda^{-1}\bar{B}_l\bar{B}_m\partial_{lm}^2 w
=0,\\{\rm div}w=0,\end{array}\right.\eqno(1.40)$$with the corresponding jump conditions\setcounter{equation}{40}
\begin{eqnarray}&&\llbracket w\rrbracket=0,\quad \llbracket -\mu (Dw+Dw^T) + \tilde{q} I\rrbracket e_d =\lambda^{-1}g[\rho]w_de_d+\lambda^{-1}\bar{B}_d\bar{B}_l\llbracket\partial_{l}w\rrbracket,
\end{eqnarray}and the boundary conditions
\begin{equation}w_-(t,x',-1)=w_+(t,x',1)=0.
\end{equation}

We make the further structural
assumption that the $x'$ dependence of $w,\ \tilde{q}$ is given as a Fourier
mode ${\rm e}^{ix'\cdot\xi}$ for $\xi\in \mathbb{R}^{d-1}$. Together with the
growing mode ansatz, this constitutes a "normal mode"  ansatz. At the rest of this subsection, we shall write down the analysis for the three dimension $d=3$,
the case for $d=2$ can be tracked readily. We define the new
unknowns $\varphi,\theta,\psi,\pi: (-1,1)\rightarrow \mathbb{R}$
by
\begin{equation} w_1(x)=-i\varphi( x_3){\rm e}^{ix'\cdot\xi},\ w_2(x)=-i\theta( x_3){\rm e}^{ix'\cdot\xi},\ w_3(x)= \psi( x_3){\rm e}^{ix'\cdot\xi},\ \tilde{q}(x)=\pi(x_3){\rm e}^{ix'\cdot\xi}. \end{equation}
To write down the equations for $\varphi,\theta,\psi,\pi$, we first notice that, denoting $'=d/dx_3$,
\begin{equation} (Dw+Dw^T)e_3=(i(\xi_1\psi-\varphi'),\ i(\xi_2\psi-\theta'),\ 2\psi')^T. \end{equation}
The equations will be quite different for the cases  that $\bar{B}$ is vertical and horizontal.

We first treat  the vertical case that $\bar{B}=(0,0,B)$. In this case, we  deduce from the equations (1.40) that  for each fixed nonzero spatial frequency $\xi=(\xi_1,\xi_2)$, and for $\varphi,\theta,\psi,\pi$ and
$\lambda$  we arrive at the following system of ODEs
$$\left\{\begin{array}{ll} \lambda^2\rho\varphi-\lambda\xi_1\pi+\mu\lambda(|\xi|^2\varphi-\varphi'')-|B|^2\varphi''=0,
\\\lambda^2\rho\theta-\lambda\xi_2\pi+\mu\lambda(|\xi|^2\theta-\theta'')-|B|^2\theta''=0,
\\\lambda^2\rho\psi+\lambda\pi'+\mu\lambda(|\xi|^2\psi-\psi'')-|B|^2\psi''=0,
\\ \xi_1\varphi+\xi_2\theta+\psi'=0,
\end{array}\right.\eqno(1.45)$$
along with the jump conditions
$$\left\{\begin{array}{ll}
\llbracket\varphi\rrbracket=\llbracket\theta\rrbracket=\llbracket\psi\rrbracket=0,
\\\llbracket\mu\lambda(\xi_1\psi-\varphi')+|B|^2\varphi'\rrbracket=\llbracket\mu\lambda(\xi_2\psi-\theta')+|B|^2\theta'\rrbracket=0,
\\\llbracket-2\mu\lambda\psi'+\lambda\pi-|B|^2\psi'\rrbracket= g[\rho]\psi,
\end{array}\right.\eqno(1.46)$$\setcounter{equation}{46}
and the boundary conditions
\begin{equation}
\varphi(-1)=\varphi(1)=\theta(-1)=\theta(1)=\psi(-1)=\psi(1)=0.\end{equation}

Now eliminating $\pi$ from the third equation in $(1.45)$ by using the other equations, we arrive at the following fourth-order ODE for $\psi$
\begin{eqnarray}-\lambda^2\rho(|\xi|^2\psi-\psi'')=\mu\lambda(|\xi|^4\psi-2|\xi|^2\psi''+\psi'''')+|B|^2(-|\xi|^2\psi''+\psi''''),
\end{eqnarray}
along with the jump conditions
\begin{eqnarray} &&
\llbracket\psi\rrbracket=\llbracket\psi'\rrbracket=0
\\&&\llbracket\mu\lambda(|\xi|^2\psi+{\psi''})+|B|^2{\psi''}\rrbracket=0
\\
&&\llbracket\mu\lambda(\psi'''-3|\xi|^2\psi')+|B|^2{\psi'''}\rrbracket=\llbracket\lambda^2\rho{\psi'}\rrbracket+g[\rho]|\xi|^2\psi.
\end{eqnarray}
and the boundary conditions
\begin{equation}
\psi(-1)=\psi(1)=\psi'(-1)=\psi'(1)=0.\end{equation}

Consequently, for $\bar{B}$ is vertical to look for growing normal mode solutions  of the original linearized problem (1.35)--(1.37) reduces to find solutions of the ordinary differential system (1.48)--(1.52). Hence, we may say in this paper that if there exists such a pair $(\lambda, \psi)$ with $\lambda>0$ satisfying (1.48)--(1.52), then the linearized problem (1.35)--(1.37) with vertical $\bar{B}$  is unstable; otherwise, the linearized problem is stable.

Now we treat  the horizontal case that $\bar{B}=(B,0,0)$, without loss of generality. We deduce from the equations (1.40)  the following system of ODEs for $\varphi,\theta,\psi,\pi$ and
$\lambda$
$$\left\{\begin{array}{ll} \lambda^2\rho\varphi-\lambda\xi_1\pi+\mu\lambda(|\xi|^2\varphi-\varphi'')+|B|^2\xi_1^2\varphi=0,
\\\lambda^2\rho\theta-\lambda\xi_2\pi+\mu\lambda(|\xi|^2\theta-\theta'')+|B|^2\xi_1^2\theta=0,
\\\lambda^2\rho\psi+\lambda\pi'+\mu\lambda(|\xi|^2\psi-\psi'')+|B|^2\xi_1^2\psi=0,
\\ \xi_1\varphi+\xi_2\theta+\psi'=0,
\end{array}\right.\eqno(1.53)$$
along with the jump conditions
$$\left\{\begin{array}{ll}
\llbracket\varphi\rrbracket=\llbracket\theta\rrbracket=\llbracket\psi\rrbracket=0,
\\\llbracket\mu\lambda(\xi_1\psi-\varphi')\rrbracket=\llbracket\mu\lambda(\xi_2\psi-\theta')\rrbracket=0,
\\\llbracket-2\mu\lambda\psi'+\lambda\pi\rrbracket= g[\rho]\psi,
\end{array}\right.\eqno(1.54)$$\setcounter{equation}{54}
and the boundary conditions
\begin{equation}
\varphi(-1)=\varphi(1)=\theta(-1)=\theta(1)=\psi(-1)=\psi(1)=0.\end{equation}
Eliminating $\pi$ from the third equation in $(1.53)$ we obtain  the following ODE for $\psi$
\begin{eqnarray}-\lambda^2\rho(|\xi|^2\psi-\psi'')=\mu\lambda(|\xi|^4\psi-2|\xi|^2\psi''+\psi'''')+|B|^2(|\xi|^2\xi_1^2\psi-\xi_1^2\psi'')
\end{eqnarray}
along with the jump conditions
\begin{eqnarray} &&
\llbracket\psi\rrbracket=\llbracket\psi'\rrbracket=0,
\\&&\llbracket\mu\lambda(|\xi|^2\psi+{\psi''})\rrbracket=0,
\\
&&\llbracket\mu\lambda(\psi'''-3|\xi|^2\psi')\rrbracket=\llbracket\lambda^2\rho{\psi'}\rrbracket+g[\rho]|\xi|^2\psi,
\end{eqnarray}
and the boundary conditions
\begin{equation}
\psi(-1)=\psi(1)=\psi'(-1)=\psi'(1)=0.\end{equation}
\section{Main results}

We first remark that if $\mu=\bar{B}=0$, $i.e.$ for the two-phase free boundary problem of the incompressible Euler equations,
the equation for $\psi$ is replaced by
\begin{eqnarray} \rho|\xi|^2\psi-\rho\psi''=0,
\end{eqnarray}
along with the modified boundary conditions
\begin{eqnarray} &&\psi(-1)=\psi(1)=0,\quad
\llbracket\psi\rrbracket=0,\quad \llbracket\rho{\psi'}\rrbracket=-\lambda^{-2}g[\rho]|\xi|^2\psi.
\end{eqnarray}
We directly obtain from (2.1)--(2.2) that
\begin{eqnarray} \int_{-1}^1\rho(|\xi|^2\psi+|\psi'|^2)\,dx_3=\lambda^{-2}g[\rho]|\xi|^2\psi^2(0).
\end{eqnarray}This immediately implies that if $[\rho]>0$, then nontrivial solutions $\psi(\xi,x_3)$ with $\lambda(\xi)>0$  can be found.
Moreover,  $\lambda(\xi)$ can be chosen as $\lambda(\xi)\ge C|\xi|$ which implies $\lambda(\xi)$ growing arbitrarily as $\xi\rightarrow\infty$. This leads to the classical ill-posedness of the Rayleigh-Taylor problem for the  Euler equations without surface tension, see \cite{E,GT1}.

Notice that in the absence of viscosity (of course with the boundary conditions modified accordingly) and for any fixed spatial frequency $\xi\neq0$, (1.48)--(1.52) (or (1.56)--(1.60)) can be viewed as an eigenvalue problem with eigenvalue $-\lambda^2$. It is not hard to check that such
a problem has a natural variational structure that allows for the construction of
solutions via the direct methods, see below. However, the presence of the viscosity  results the
appearance of $\lambda$ both quadratically and linearly which
destroys the variational structure of these eigenvalue problems. To restore the ability to use variational methods,
 the
authors in \cite{GT2} developed a quite general and robust method. More precisely, they  artificially remove the linear terms
$\mu\lambda$ by  the modified viscosities
$\tilde{\mu}=s\mu$
where $s>0$ is an arbitrary parameter. After establishing the
existence of the solutions to the modified problem, by proving that there
is a fixed point $\lambda=s$, then the corresponding solution is a
solution to the original problem. We use the same trick in this paper to construct solutions. However, the purpose of this paper is to study in detail the stabilized effect of the magnetic field in the Rayleigh-Taylor instability, to characterize the critical magnetic number and
the critical frequency when the magnetic number under the critical value.

We assume in this paper that $\mu, B$ are not zero. We define the
critical magnetic number $|B|_c$ through the following variational
problem
\begin{equation}
|B|_c^2:=\sup_{\psi\in
H_0^1((-1,1))}\frac{g[\rho]\psi^2(0)}{\int_{-1}^1|\psi'|^2\,dy},\end{equation}
which only depends on $g,\ [\rho]$. For $|B|<|B|_c$, we define the critical frequency
$|\xi|_{vc}^B$ for the vertical $\bar{B}$ by the variational problem
\begin{equation}
(|\xi|^B_{vc})^2:=\inf_{\psi\in
H_0^2((-1,1))}\frac{|B|^2\int_{-1}^1|\psi''|^2\,dy}{g[\rho]\psi^2(0)-|B|^2\int_{-1}^1|\psi'|^2\,dy},\end{equation}
and define the critical frequency
$|\xi|_{hc}^B$ for the horizontal $\bar{B}$ by the variational problem
\begin{equation}
(|\xi|^B_{hc})^2:=\sup_{\psi\in
H_0^1((-1,1))}\frac{g[\rho]\psi^2(0)-|B|^2\int_{-1}^1|\psi'|^2\,dy}{|B|^2\int_{-1}^1|\psi|^2\,dy}.\end{equation}
The
superscript emphasizes the dependence of the critical frequencies (2.5)--(2.6) on
$|B|$, these two only depend on $g,\ [\rho]$ and $|B|$.  We will show in Lemma 3.2 that the extremums (2.4)--(2.6)
are achieved by the direct method.

Now we may first state out main results of the paper for the two-dimensional case, the results for three-dimensional case will be discussed later. The first one is concerned with the problem (1.48)--(1.52) for the vertical magnetic number.
\begin{theorem} Let $\bar{B}=(0,B)$ be vertical, then we have:

$(i)$ For any fixed $|B|\ge|B|_c$ and any $\xi\in \mathbb{R}$ there
is no nontrivial solution $\psi $ with $\lambda >0$ to the problem
(1.48)--(1.52);

$(ii)$ For any fixed $|B|<|B|_c$, then for $\xi\in \mathbb{R}$ so
that $ |\xi|\le |\xi|_{vc}^B$ there is no nontrivial solution $\psi $
with $\lambda >0$ to the problem
(1.48)--(1.52);

$(iii)$ For any fixed $|B|<|B|_c$, then for $\xi\in \mathbb{R}$ so
that $|\xi|> |\xi|_{vc}^B$ there exists
$\psi=\psi(\xi,x_2)$ and $\lambda( \xi
)>0$ to the problem
(1.48)--(1.52). Moreover, $\psi,\ \lambda$ are even in $\xi$ and $\psi$ is smooth when
restricted to $(-1,0)$ or $(0,1)$ with $\psi(\xi,0)\neq 0$.
\end{theorem}

Next theorem is concerned with the problem (1.56)--(1.60) for the horizontal magnetic number.
Note that since we consider the 2D case, $\xi_1^2=|\xi|^2$ for the moment.
\begin{theorem} Let $\bar{B}=(B,0)$ be horizontal, then we have:

$(i)$ For any fixed $|B|\ge|B|_c$ and any $\xi\in \mathbb{R}$ there
is no nontrivial solution $\psi $ with $\lambda >0$ to the problem
(1.56)--(1.60).

$(ii)$ For any fixed $|B|<|B|_c$, then for $\xi\in \mathbb{R}$ so
that $ |\xi|\ge |\xi|_{hc}^B$ there is no nontrivial solution $\psi $
with $\lambda >0$ to the problem
(1.56)--(1.60).

$(iii)$ For any fixed $|B|<|B|_c$, then for $\xi\in \mathbb{R}$ so
that $|\xi|< |\xi|_{hc}^B$ there exists
$\psi=\psi(\xi,x_2)$ and $\lambda( \xi
)>0$ to the problem
(1.56)--(1.60). Moreover, $\psi,\ \lambda$ are even in $\xi$ and $\psi$ is smooth when
restricted to $(-1,0)$ or $(0,1)$ with $\psi(\xi,0)\neq 0$.
\end{theorem}

As we will see, the conclusions of Theorem 2.1 and Theorem 2.2 also hold for the inviscid fluids ($i.e.,\ \mu=0$). They demonstrate that the stabilized effect of the magnetic field is much more remarkable than  that of the surface tension $\sigma$ which only stabilizes the frequency interval $(0, \sqrt{\frac{g[\rho]}{\sigma}})$ for any $\sigma>0$, see \cite{C,GT2}.

It is easy to imagine the corresponding results for the 3D case. Indeed, for the vertical magnetic number $\bar{B}=(0,0, B)$, the problem
(1.48)--(1.52) is the same and hence Theorem 2.1 holds. However, for the horizontal $\bar{B}=(B,0,0)$, if we consider the frequency $\xi=(\xi_1,\xi_2)$
with $\xi_1=0$. Then (1.56)--(1.60) becomes (1.48)--(1.52) with $B=0$, there is no effect of the magnetic number and hence by Theorem 2.1 $(iii)$ the system is unstable. In other words, the horizontal magnetic number only stabilizes the frequencies along its own direction.

It is important to know the behavior of the eigenvalue $\lambda(\xi)$ obtained in $(iii)$ of both Theorem 2.1 and Theorem 2.2 within their respective definitional intervals.
\begin{theorem} Let $0<|B|<|B|_c$, then we have

(i) Let $\lambda(\xi)$ be that of $(iii)$ in Theorem 2.1, then
$\lambda:(|\xi|_{vc}^B,\infty)\rightarrow(0,\infty)$ is continuous and
\begin{equation}\lim_{|\xi|\rightarrow|\xi|_{vc}^B}\lambda(|\xi|)=0\hbox{ and }\sup_{|\xi|_{vc}^B<|\xi|<\infty}\lambda(|\xi|)\le\frac{2\sqrt{g[\rho]}}{|B|\sqrt[4]{\rho_+}} .\end{equation}

(ii) Let $\lambda(\xi)$ be that of $(iii)$ in Theorem 2.2, then
$\lambda:(0,|\xi|_{hc}^B)\rightarrow(0,\infty)$ is continuous and
\begin{equation}\lim_{|\xi|\rightarrow0}\lambda(|\xi|)=\lim_{|\xi|\rightarrow|\xi|_{hc}^B}\lambda(|\xi|)=0.\end{equation}
\end{theorem}

It is worth to point out that the proof of Theorem 2.3 is independent of the viscosity $\mu$ and Theorem 2.3 allows us to define the fastest growth rate by, for vertical and horizontal magnetic number, respectively,
\begin{equation}\Lambda_v=\sup_{|\xi|_{vc}^B<|\xi|<\infty}\lambda(|\xi|)<\infty\hbox{ and }\Lambda_h=\sup_{0<|\xi|<|\xi|_{hc}^B}\lambda(|\xi|)<\infty.\end{equation}
This means that the growing mode solutions constructed in Theorem 2.1 and Theorem 2.2 (except for the horizontal $\bar{B}$ in 3D case) do not grow arbitrarily, which may imply that the magnetic field prevent the ill-posedness for the inviscid Rayleigh-Taylor problem, see $\cite{E,GT1}$.

With $\psi$ constructed in  $(iii)$ of Theorem 2.1 (resp. $(iii)$ of Theorem 2.2), one can use the system of ODEs (1.45)--(1.47) (resp. (1.53)--(1.55)) to obtain a growing normal mode solution to the linearized problem (1.35)--(1.37) for $\bar{B}$ is vertical (resp. $\bar{B}$ is horizontal). Although ${\rm e}^{ix'\xi}\notin L^2(\Omega)$, as in \cite{GT1,GT2}, we can resort to a Fourier synthesis of such solution to construct
 solutions to (1.35)--(1.37) which grow in the piecewise Sobolev space $H^k,\ \forall k\ge 1$. The growth rate of such growing solutions is not exactly the fastest growth rate ${\rm e}^{t\Lambda_v}$ (resp. ${\rm e}^{t\Lambda_h}$),  but arbitrarily close to this rate. On the other hand, based on the proof of Theorem 2.1--2.3, we can estimate the growth in time of arbitrary solutions to (1.35)--(1.37) in terms of ${\rm e}^{t\Lambda_v}$ (resp. ${\rm e}^{t\Lambda_h}$). The statement and the proof of these results is very similar to that of \cite[Theorem 2.4--2.5]{GT2}, we omit them and let the interested readers refer to \cite{GT2} for the details.

Notice that in Theorem 2.1--2.2, when $|B|\ge |B|_c$ no growing mode solutions can be constructed by our argument. This should suggest that the linear system (1.35)--(1.37) is stable when $|B|\ge |B|_c$. In fact, when $|B|>|B|_c$ we can prove the following stability estimates:
\begin{theorem} Let $|B|>|B|_c$ be fixed, then

$(i)$ for $\bar{B}=(0,B)$, we have\begin{eqnarray}&&
\|\partial_tv(t)\|_{L^2}^2+\|v(t)\|_{L^2}^2+\|\partial_2v(t)\|_{L^2}^2+\int_0^t\|\partial_t v(s)\|_{H^1}^2\,ds
\nonumber\\&&\quad\le C( \|\partial_tv(0)\|_{L^2}^2+\|v(0)\|_{L^2}^2+\|\partial_2v(0)\|_{L^2}^2);
\end{eqnarray}
if in addition, $\eta(0,x_1,-1)=\eta(0,x_1,1)=0$ initially, then we have
\begin{eqnarray}&&
\|v(t)\|_{L^2}^2+\|\eta(t)\|_{L^2}^2+\|\partial_2\eta(t)\|_{L^2}^2+\int_0^t\|v(s)\|_{H^1}^2\,ds
\nonumber\\&&\quad\le C( \|v(0)\|_{L^2}^2+\|\eta(0)\|_{L^2}^2+\|\partial_2\eta(0)\|_{L^2}^2).
\end{eqnarray}

$(ii)$ for $\bar{B}=(B,0)$, we have
\begin{eqnarray}&&
\|\partial_tv(t)\|_{L^2}^2+\|\partial_1v_1(t)\|_{L^2}^2+\|v_2(t)\|_{H^1}^2+\int_0^t\|\partial_tv(s)\|_{H^1}^2\,ds
\nonumber\\&&\quad\le C( \|\partial_tv(0)\|_{L^2}^2+\|\partial_1v_1(0)\|_{L^2}^2+\|v_2(0)\|_{H^1}^2);
\end{eqnarray}
if in addition, $(\partial_1\eta_1+\partial_2\eta_2)(0)=0$ initially, then we have
\begin{eqnarray}&&
\|v(t)\|_{L^2}^2+\|\partial_1\eta_1(t)\|_{L^2}^2+\|\eta_2(t)\|_{H^1}^2+\int_0^t\|v(s)\|_{H^1}^2\,ds
\nonumber\\&&\quad\le C( \|v(0)\|_{L^2}^2+\|\partial_1\eta_1(0)\|_{L^2}^2+\|\eta_2(0)\|_{H^1}^2).
\end{eqnarray}

\end{theorem}

As before, the assertion $(i)$ also holds for 3D case. The proof of Theorem 2.4 follows by combining the standard energy estimates and the definition (2.4) of the critical magnetic number $|B|_c$. It is understood that the initial value $\partial_tv(0)$ are given in terms of the initial data through the equations $(1.35)_2$.
 It is obvious that the estimates (2.10)--(2.12) also hold for $\partial_t^i\partial_1^j,\ \forall i,j\ge 1$ since they satisfy the same system  (1.35)--(1.37). It in particular implies that the exponential growth of solutions is impossible.

At last, we emphasize that the success of deriving the critical magnetic number (and hence critical frequency) essentially depends  on the fact that we consider the problem within
a slab. The critical magnetic number does not hold for the whole space problem and indeed in this case for any $\bar{B}$ there always exists growing normal mode solutions, see \cite{C}.

The instability of certain steady states of different densities which occurs under an acceleration
of the fluid system in the direction toward the denser fluid is well known as the Rayleigh-Taylor instability since the works \cite{R,T}. It attracts huge attention of the researches both physically and mathematically, we refer to book \cite{C} and the report \cite{K} for the detailed references.
Our construction of the growing normal mode solutions by variational method is inspired by the works \cite{GH,GT1,GT2,H} which deal with the nonconstant density profile.

The rest of the paper is devoted to prove our theorems.

\section{Proof of Theorems}
\subsection{Proof of Theorem 2.1}
We treat the problem (1.48)--(1.52) to prove Theorem 2.1 in this subsection. As quoted in Section 2, to restore the ability to use variational methods we remove the linear dependence on $\lambda$ by defining  the modified viscosities $\tilde{\mu}=s\mu$,
where $s>0$ is an arbitrary parameter. Then we obtain a family $(s>0)$ of
modified problems
\begin{eqnarray}-\lambda^2\rho(|\xi|^2\psi-\psi'')=s\mu(|\xi|^4\psi-2|\xi|^2\psi''+\psi'''')+|B|^2(-|\xi|^2\psi''+\psi''''),
\end{eqnarray}
along with the jump conditions
\begin{eqnarray} &&
\llbracket\psi\rrbracket=\llbracket\psi'\rrbracket=0,
\\&&\llbracket s\mu(|\xi|^2\psi+{\psi''})+|B|^2{\psi''}\rrbracket=0,
\\
&&\llbracket s\mu(\psi'''-3|\xi|^2\psi')+|B|^2{\psi'''}\rrbracket=\llbracket\lambda^2\rho{\psi'}\rrbracket+g[\rho]|\xi|^2\psi.
\end{eqnarray}
and the boundary conditions
\begin{equation}
\psi(-1)=\psi(1)=\psi'(-1)=\psi'(1)=0.\end{equation}

Notice that for any fixed $s>0$ and $\xi$, (3.1)--(3.5) is a standard eigenvalue problem for $-\lambda^2$. It allows us to use the variational methods to construct
solutions. We define the energies
\begin{eqnarray}&&E(\psi)
=\frac{1}{2}\int_{-1}^1s\mu(4|\xi|^2|\psi'|^2+||\xi|^2\psi+\psi''|^2)+|B|^2(|\xi|^2|\psi'|^2+|\psi''|^2)\,dx_2\nonumber
\\&&\qquad\qquad-\frac{1}{2}|\xi|^2g[\rho]|\psi(0)|^2,
\\&&J(\psi)
=\frac{1}{2}\int_{-1}^1 \rho(|\xi|^2|\psi|^2+|\psi'|^2)\,dx_2,\end{eqnarray}
  which are both
well-defined on the space $H_0^2((-1,1))$, the subset of $H^2((-1,1))$ satisfying (3.5). We define the set
\begin{eqnarray}\mathcal{A}=\{ \psi \in
H_0^2((-1,1)) \ |\
J(\psi)=1\},\end{eqnarray} then we want to find the smallest
$-\lambda^2$ by minimizing
\begin{equation}-\lambda^2(|\xi|)=\alpha(|\xi|):=\inf_{\psi\in
\mathcal{A}}E(\psi).\end{equation}

The first thing is to show that a minimizer of (3.9) exists and that the minimizer satisfies Euler-Langrange
equations equivalent to (3.1)--(3.5).
\begin{Proposition} (i)For any fixed $\xi\neq0$ and $s>0$, $E$ achieves its infimum on $\mathcal{A}$.

(ii)Let
$\psi$ be a minimizer and $-\lambda^2=\alpha:=E(\psi)$, then the pair $\psi,\lambda^2$
satisfies (3.1) along with the jump and boundary conditions
(3.2)--(3.5). Moreover, $\psi$ is smooth when restricted to
$(-1,0)$ or $(0,1)$.
\end{Proposition}
\hspace{-18pt}{\bf Proof.} To prove $(i)$ we first note that for any $\psi\in \mathcal{A}$
\begin{eqnarray} E(\psi)\ge -\frac{1}{2}|\xi|^2[\rho]|\psi(0)|^2\ge -\frac{1}{2}|\xi|^2[\rho]\int_{0}^1 |\psi'|^2\,dx_2\ge -({\rho_+})^{-1}[\rho]{|\xi|^2}.\end{eqnarray}
Hence $E$ is bounded below on $\mathcal{A}$. Let $\psi_n\in\mathcal{A}$ be a minimizing sequence, then  $E(\psi_n)$ is  bounded. This together with (3.6), and (3.10) again imply that $\psi_n$ is bounded in $H_0^2((-1,1))$. So there exists $\psi\in H_0^2((-1,1))$ such that $\psi_n\rightarrow\psi$
weakly in $\psi\in H_0^2((-1,1))$ and strongly in $C^1((-1,1))$. $\psi\in \mathcal{A}$ follows from the strong convergence. Moreover, by the lower semi-continuity and the strong convergence, we have
\begin{eqnarray}E(\psi)\le \liminf_{n\rightarrow\infty}E(\psi_n)=\inf_\mathcal{A} E. \end{eqnarray}

Now we prove $(ii)$. We notice that since $E$ and $J$ are homogeneous of degree 2, (3.9) is equivalent to
\begin{equation}-\lambda^2=\inf_{\psi\in
H_0^2((-1,1))}\frac{E(\psi)}{J(\psi)}.\end{equation}
For $\tau\in \mathbb{R}$ and any $\psi_0\in H_0^2((-1,1))$, define $\psi(\tau)=\psi+\tau\psi_0$, then (3.12) implies
\begin{equation}E(\psi(\tau))+\lambda^2J(\psi(\tau))\ge 0.\end{equation}
If we set $I(\tau)=E(\psi(\tau))+\lambda^2J(\psi(\tau))$, then we have $I(\tau)\ge 0$ for all $\tau\in \mathbb{R}$ and $I(0)=0$. This implies $I'(0)=0$. By the expressions (3.6)--(3.7), direct computation leads to
\begin{eqnarray}&&
\int_{-1}^1s\mu(4|\xi|^2\psi'\psi_0'+ (|\xi|^2\psi+\psi'')(|\xi|^2\psi_0+\psi_0''))+|B|^2(|\xi|^2\psi'\psi_0'+\psi''\psi_0'')\,dx_2\nonumber
\\&&\quad=|\xi|^2[\rho]\psi(0)\psi_0(0)-\lambda^2\int_{-1}^1 \rho(|\xi|^2\psi\psi_0+\psi'\psi_0')\,dx_2.\end{eqnarray}

By further assuming $\psi_0$ is compactly supported in either $(-1,0)$ or $(0,1)$, we find that $\psi$ satisfies the equation (3.1) in a weak
sense on $(-1,0)$ and $(0,1)$. A standard bootstrap argument then shows that $\psi$ is in $H^k((-1,0))$ (resp. $H^k((0,1))$) for all $k\ge 0$
when restricted to $(-1,0)$ (resp. $(0,1)$), and hence is smooth when restricted to either interval. Since $\psi\in H_0^2((-1,1))$, the jump conditions (3.2) and the boundary conditions (3.5) follow trivially. It remains to show that the jump conditions (3.3) and (3.4) hold. For this we take $\psi_0\in C_c^\infty((-1,1))$ in (3.14). Integrating (3.14) by parts the terms of $\psi_0''$, we obtain
\begin{eqnarray}&&
\int_{-1}^1s\mu(|\xi|^4\psi \psi_0+3|\xi|^2\psi'\psi_0'+|\xi|^2\psi''\psi_0- \psi'''\psi_0')+|B|^2(|\xi|^2\psi'\psi_0'-\psi'''\psi_0')
\,dx_2
\\&&\quad=\llbracket s\mu(|\xi|^2\psi+\psi'')+|B|^2\psi''\rrbracket\psi_0'(0)+|\xi|^2[\rho]\psi(0)\psi_0(0)-\lambda^2\int_{-1}^1 \rho(|\xi|^2\psi\psi_0+\psi'\psi_0')\,dx_2.\nonumber\end{eqnarray}
Integrating further (3.15) by parts the  terms of $\psi_0'$,
using the fact that $\psi$ solves (3.1) on $(-1,0)$ and $(0,1)$ and the jump conditions (3.2), we find that
\begin{eqnarray}&&-\llbracket s\mu(|\xi|^2\psi+\psi'')+|B|^2\psi''\rrbracket\psi_0'(0)+
\llbracket s\mu(\psi'''-3|\xi|^2\psi')+|B|^2{\psi'''}\rrbracket\psi_0(0)\nonumber
\\&&\quad=\left(\llbracket\lambda^2\rho{\psi'}\rrbracket+g[\rho]|\xi|^2\psi(0)\right)\psi_0(0).\end{eqnarray}
Since $\psi_0$ may be chosen arbitrarily, (3.16) implies the jump conditions (3.3)--(3.4) and we conclude our proof.\hfill$\Box$\smallskip\smallskip

Now we come to the heart part of this paper to clarify the sign of
the infimum obtained in Proposition 3.1. It is crucial to represent the energy
$E(\psi)$ in the following form
\begin{equation}E(\psi)=|\xi|^2E_0(\psi)+sE_1(\psi),\end{equation}
where
\begin{eqnarray}&&E_0( \psi)=\frac{1}{2}\int_{-1}^1|B|^2\left(|\psi'|^2+\frac{|\psi''|^2}{|\xi|^2}\right)\,dx_2-\frac{1}{2}g[\rho]|\psi(0)|^2,
\\&&E_1( \psi)
=\frac{1}{2}\int_{-1}^1\mu(4|\xi|^2|\psi'|^2+||\xi|^2\psi+\psi''|^2) \,dx_2.
\end{eqnarray}
Since the parameter $s$ is positive but can be made to be small arbitrarily, the key
point is to clarify the sign of the energy $E_0(\psi)$. Note that if without magnetic field
this energy is non-positive, hence the pure fluid Rayleigh-Taylor problem is unstable.
However, the presence of the magnetic field makes it have possibility to be positive.
It
is characterized by the critical magnetic value $|B|_c$ and the
critical frequency $|\xi|_{vc}^B$ (or $|\xi|_{hc}^B$ for the horizontal case), we recall their definitions (2.4)--(2.6).

\begin{lemma} $(i)$ The supremum  (2.4)  is achieved;

$(ii)$ For any fixed $B$ so that $|B|< |B|_c$, the infimum  (2.5)
is achieved. Moreover, $|\xi|_{vc}^B$ is continuous and strictly
decreasing as a function of $|B|$ within $0<|B|< |B|_c$ and we have
\begin{equation}|\xi|_{vc}^B\rightarrow0\ \hbox{ as }\ |B|\rightarrow 0,\  \hbox{ and }\ |\xi|_{vc}^B\rightarrow\infty\ \hbox{
as }\  |B|\rightarrow|B|_c.
\end{equation}

$(iii)$ For any fixed $B$ so that $|B|< |B|_c$, the supremum  (2.6)
is achieved. Moreover, $|\xi|_{vc}^B$ is continuous and strictly
decreasing as a function of $|B|$ within $0<|B|< |B|_c$ and we have
\begin{equation}|\xi|_{hc}^B\rightarrow\infty\ \hbox{ as }\ |B|\rightarrow 0,\  \hbox{ and }\ |\xi|_{hc}^B\rightarrow0\ \hbox{
as }\  |B|\rightarrow|B|_c.
\end{equation}
\end{lemma}
\hspace{-18pt}{\bf Proof.} To prove $(i)$, we express the supremum
(2.4) in an equivalent form
\begin{equation}
\frac{1}{|B|_c^2}:=\inf\int_{-1}^1|\psi'|^2\,dx_2 ,\end{equation}
where $ \psi\in H_0^1((-1,1))$ satisfying the constraint
\begin{equation}
 {g[\rho]}\psi^2(0) ={1}.
\end{equation}
It is obvious that the integral in (3.22) is bounded below. Let
$\psi_n $ be a minimizing sequence, then $\psi_n$ is bounded in
$H_0^1$ and hence $\psi_n\rightarrow\psi$ weakly in $H_0^1$, up to
the extraction of a subsequence if necessary. Owing to the compact
embedding $H_0^1\hookrightarrow C^0$ we have
$\psi_n(0)\rightarrow\psi(0)$. By these convergences and the weak
semi-continuity of the integral in (3.22), we obtain that $\psi$ is
a minimizer of (3.22) and this proves $(i)$.

 Now we prove $(ii)$. To prove that the infimum (2.5) is achieved,
 we first
rewrite (2.5) equivalently
\begin{equation}
(|\xi|^B_{vc})^2:=\inf_{\psi\in \mathcal{C}}P(\psi),\end{equation}
where
\begin{eqnarray}&&P(\psi)=|B|^2\int_{-1}^1|\psi''|^2\,dx_2,
\\&&
\mathcal{C}=\left\{\psi\in  H_0^2((-1,1))\ \left|\
{{g[\rho]\psi^2(0)-|B|^2\int_{-1}^1|\psi'|^2\,dx_2}}\right.=1\right\}.\end{eqnarray}
Let $\psi_n\in
\mathcal{C}$ be a minimizing sequence, we have from (3.24) that $\psi_n''$ is bounded in $L^2$. This with Poincar\'e's inequality
implies that $\psi_n$ is bounded in $H_0^2$. So we have $\psi_n\rightarrow\psi$
weakly in $H_0^2$, up to the extraction of a subsequence if
necessary, and owing to the compact embeddings $H_0^2\hookrightarrow C^1$ we have $\psi_n\rightarrow\psi$
strongly in $H^1$ and $\psi_n(0)\rightarrow\psi(0)$ as well. These
convergences and the weak semi-continuity yield
\begin{equation}P(\psi)\le  \limsup_{n\rightarrow\infty}P(\psi_n)=\inf_\mathcal{C}P.\end{equation}
$\psi\in \mathcal{C}$ follows from the strong convergences and
hence $\psi$ is a maximizer of (3.22). On the other hand,
by the definitions (2.4) and (2.5), it is easy to have the
conclusions of the behavior of $|\xi|_{vc}^B$ as $B$ varies within
$0<|B|<|B|_c$.

 At last, we prove $(iii)$. We first rewrite the supremum (2.6) equivalently
\begin{equation}
(|\xi|^B_{hc})^2:=\sup_{\psi\in \mathcal{F}}Q(\psi) ,\end{equation}
where
\begin{eqnarray}&&
Q(\psi)=
{g[\rho]\psi^2(0)-|B|^2\int_{-1}^1|\psi'|^2\,dx_2} ,
\\&&\mathcal{F}=\left\{\psi\in H_0^1((-1,1))\ \left|\
{|B|^2\int_{-1}^1|\psi|^2\,dx_2}\right.=1\right\}.\end{eqnarray} We
first show that $Q$ is bounded above on $\mathcal{F}$. Indeed, for
any $\psi\in \mathcal{F}$, we have
\begin{eqnarray}&Q(\psi)&=
{2g[\rho]\psi^2(0)-|B|^2\int_{-1}^1|\psi'|^2\,dx_2}\nonumber
\\&&\le
C\int_0^l|\partial_{x_2}(\psi^2)|\,dx_2-|B|^2\int_{-1}^1|\psi'|^2\,dx_2\nonumber
\\&&\le
C\int_{-1}^1| \psi
|^2\,dx_2-\frac{|B|^2}{2}\int_{-1}^1|\psi'|^2\,dx_2\nonumber
\\&&\le
C\int_{-1}^1| \psi |^2\,dx_2\le C.\end{eqnarray} Let $\psi_n\in
\mathcal{F}$ be a maximizing sequence, the second inequality in
(3.31) implies that
\begin{equation}\frac{|B|^2}{2}\int_{-1}^1|\psi_n'|^2\,dx_2\le  C\int_{-1}^1|
\psi_n |^2\,dx_2-Q(\psi_n)\le C.\end{equation} Then $\psi_n$ is bounded in $H_0^1$, so we have $\psi_n\rightarrow\psi$
weakly in $H_0^1$, up to the extraction of a subsequence if
necessary, and owing to the compact embeddings $H_0^1\hookrightarrow
L^2$ and $H_0^1\hookrightarrow C^0$ we have $\psi_n\rightarrow\psi$
strongly in $L^2$ and $\psi_n(0)\rightarrow\psi(0)$ as well. These
convergences and the weak semi-continuity of the integral in $Q$, we
have that
\begin{equation}Q(\psi)\ge  \limsup_{n\rightarrow\infty}Q(\psi_n)=\sup_\mathcal{F}Q.\end{equation}
$\psi\in \mathcal{F}$ follows from the strong $L^2$ convergence and
hence $\psi$ is a maximizer of (3.28).
Finally, by the definitions (2.4) and (2.6), we conclude the behavior of $|\xi|_{vc}^B$ as $B$ varies within
$0<|B|<|B|_c$. The proof of Lemma 3.2 is completed.\hfill$\Box$\smallskip\smallskip

Now we can clarify the sign of the energy $E_0(\psi)$.

\begin{lemma} We have the following four assertions:

$(i)$ If $|B|\ge |B|_c$, then for any $\psi$ we have
$E_0(\psi)\ge0$. Moreover,
\begin{equation}E_0(\psi)\ge\frac{1}{2}\int_{-1}^1
(|B|^2-|B|_c^2)
|\psi'|^2+\frac{|B|^2|\psi''|^2}{|\xi|^2}\,dx_2.\end{equation}

$(ii)$ Fixed $|B|< |B|_c$ and $|\xi|\le |\xi|_{vc}^B$, then for any
$\psi$  we have $E_0(\psi)\ge 0$.

$(iii)$ Fixed $|B|< |B|_c$ and $|\xi|>|\xi|_{vc}^B$, then there exists
$\psi$ such that $E_0(\psi)<0$.

$(iv)$ If there is $\psi$ such that $E_0(\psi)<0$,
then $|B|< |B|_c,\ |\xi|<|\xi|_{vc}^B$ and $\psi(0)\neq 0$.
\end{lemma} \hspace{-18pt}{\bf Proof.} The first three assertions follows easily by the definitions (2.4) and (2.5). Indeed, for $(i)$, we let $B$ be so that $|B|\ge|B|_c$. By the
definition (2.4), we have
\begin{equation} |B|_c^2 \int_{-1}^1 |\psi'|^2 \,dx_2\ge g[\rho] \psi^2(0),\hbox{ for any } \psi\in H_0^1 .\end{equation}
Since the remaining term in $E_0$ is nonnegative, we verify
the assertion $(i)$.
To prove $(ii)$ and $(iii)$, we fix any $B$ so that $|B|< |B|_c$.
Again similarly by the definition (2.5), we have that if
$|\xi|\le|\xi|^B_{vc}$ then
\begin{equation}
\int_{-1}^1|B|^2(|\xi|^2|\psi'|^2+{|\psi''|^2})\,dx_2\ge g|\xi|^2[\rho]|\psi(0)|^2.
\end{equation}
This proves the assertion $(ii)$. Now if $|\xi|>
|\xi|^B_{vc}$, by the definition (2.6) there exists $\psi\in H_0^1$ such
that
\begin{equation}
\int_{-1}^1|B|^2(|\xi|^2|\psi'|^2+{|\psi''|^2})\,dx_2< g|\xi|^2[\rho]|\psi(0)|^2.
\end{equation}
This prove $(iii)$.

It remains to prove $(iv)$. If $E_0(\psi)<0$, then the
assertions $|B|< |B|_c$ and $|\xi|>|\xi|_{vc}^B$ follow from $(i)$ and
$(ii)$. Since the integrals in (3.18) are all non-negative, then we
must have $\psi(0)\neq 0$. This proves $(iv)$ and we conclude our
lemma.\hfill$\Box$\smallskip\smallskip

By Lemma 3.3 we are able to show the sign of the infimum of $E$ over $\mathcal{A}$. We write $\alpha=\alpha(s)$ to emphasize the dependence on $s\in (0,\infty)$.
\begin{lemma}
$(i)$ If either $|B|\ge |B|_c$ or $|B|<|B|_c$  with
$|\xi|\le|\xi|_{vc}^B$, then $\alpha(s)\ge 0$ for any $s\ge 0$.

$(ii)$ If both $|B|< |B|_c$ and $|\xi|>|\xi|_{vc}^B$ are satisfied,
then there exists $s_0>0$ depending on the quantities
$\rho_\pm, \mu_\pm,g, |B|,|\xi|$ so that
for $s\le s_0$ it holds that $\alpha(s)<0$.
\end{lemma}
\hspace{-18pt}{\bf Proof.} Since both $E$ and $J$ are homogeneous of
degree $2$ and $J$ is positive definite, we may reduce to clarify
the sign of the energy $E(\psi)$.

To prove $(i)$, observe that in this case we deduce from Lemma 3.3
that $E_0(\psi)\ge0$ for any $\psi\in \mathcal{A}$. Combining
this with the fact $E_1(\psi)\ge 0$, we have
$E(\psi)\ge0$ for any $s\ge0$. Hence taking the infimum we
have $\mu(s)\ge 0$. This proves $(i)$.

It remains to prove $(ii)$. In this case we know from Lemma 3.3 that
there exists $ \widetilde{\psi}$ such that
$E_0(\widetilde{\psi})<0$. Obviously, we have
\begin{equation}E
(\widetilde{\psi})=|\xi|^2E_0(\widetilde{\psi})+sE_1(\widetilde{\psi})\le
E_0(\widetilde{\psi})+sC\end{equation} for a
constant $C$ depending on
$\rho_\pm, \mu_\pm,g, |B|,|\xi|$.  Then
there exists $s_0>0$ depending on  these parameters such that for
$s\le s_0$ it holds that
$E(\widetilde{\psi})<0$. Hence the infimum
$\alpha(s)<0$ for $s \le s_0$. This proves $(ii)$ and the proof of
Lemma 3.4 is completed.\hfill$\Box$\smallskip\smallskip

Immediately, Lemma 3.4 proves the assertions $(i)$ and $(ii)$ of
Theorem 2.1 by contradiction. To finish the proof of Theorem 2.1, it suffices to prove the assertion
$(iii)$. So fixed $|B|<|B|_c$ and then $|\xi|>|\xi|_{vc}^B$. We want
to show that there is a fixed point such that $\lambda=s$. To this end, we first study the behavior of $\alpha(s)$ as a function of $s\ge 0$.
\begin{lemma}We have the following statements.

$(i)$ $\alpha(s)$ is strictly increasing;

 $(ii)$
$\mu\in C_{loc}^{0,1}((0,\infty))\cap C^0((0,\infty))$;

 $(iii)$ For any $b>|\xi|_{vc}^B$, there exist  two constants
 $C_0,C_1>0$ depending on the parameters
 $\rho_\pm, \mu_\pm,g, |B|,b$ so that
\begin{equation}\alpha(s)\le-C_0+sC_1,\quad\hbox{ for all }|\xi| \in
[b,\infty);\end{equation}

 $(iv)$ There exist  constants
 $C_2>0$ depending on
 $\rho_\pm,g$ and  $C_3>0$ depending additionally on
 $\mu_\pm, |B|, \xi$  so that
\begin{equation}\alpha(s)\ge-C_2|\xi|+sC_3.\end{equation}

\end{lemma}
\hspace{-18pt}{\bf Proof.} Recall the energy decomposition (3.17)
along with (3.18) and (3.19). It keeps the same form as in
\cite[Propostion 3.6]{GT2}, hence $(i)$ and $(ii)$ follow  in the same
way.

To prove $(iii)$, fixed $b>|\xi|_{vc}^B$, by  Lemma 3.3 $(iii)$ there
exists $ \psi_b$ such that $C_0=-E_0( \psi_b)>0$.
Then we have $E( \psi_b)\le -C_0+sC_1$ for some $C_1>0$ and
then $(iii)$ holds.

Finally, we prove $(iv)$. First observe that for any $
\psi\in \mathcal{A}$ we have
\begin{eqnarray}& -|\xi|^2g[\rho]|\psi(0)|^2
&\ge -|\xi|g[\rho]\left(\int_{0}^1|\xi||\psi|^2\,dx_2\right)^\frac{1}{2}\left(\int_{0}^1|\psi'|^2\,dx_2\right)^\frac{1}{2}\ge -C_2|\xi|.\end{eqnarray}
Since the other terms in the energy $E$ is nonnegative, we have
\begin{equation}\alpha(s)\ge-C_2|\xi|+s\inf_{\psi\in \mathcal{A}}E_1(\psi).\end{equation}
We denote by $C_3$ this positive infimum, then $(iv)$ follows and we
conclude our lemma.\hfill$\Box$\smallskip\smallskip

By Lemma 3.4 and Lemma 3.5, we then define the open set
\begin{equation}\mathcal{S}=\mu^{-1}((-\infty,0))\subset(0,\infty).\end{equation}
Note that $\mathcal{S}$ is non-empty and allows us to define
$\lambda(s)=\sqrt{-\alpha(s)}$ for $s\in\mathcal{S}$. We state the
existence of solutions to the modified problem (3.1)--(3.5) which we
have already proved.
\begin{Proposition}For each $s\in \mathcal{S}$ there is a solution $\psi=\psi_s( \xi ,x_2)$
with $\lambda=\lambda(   \xi ,s)>0$ to the problem (3.1)
along with the jump and boundary conditions (3.2)--(3.5). Moreover, $\psi_s,\ \lambda$ are even in $\xi$.
The solutions are smooth when restricted to $(-1,0)$ or $(0,1)$ with
$\psi_s(\xi,0)\neq 0$.
\end{Proposition}
\hspace{-18pt}{\bf Proof.} Let $\psi_s( \xi ,x_2)$ be constructed in Proposition 3.1. Since $s\in \mathcal{S}$, we can write
 $-\alpha(\xi,s)=\lambda^2(\xi,s)$, then $\psi_s( \xi ,x_2),\lambda(\xi ,s)$ solve the problem (3.1)--(3.5). The remaining assertions follow from Proposition 3.1 $(ii)$ and Lemma 3.3 $(iv)$. The proof of lemma is completed.\hfill$\Box$\smallskip\smallskip

Now we will make a fixed-point argument to find $s\in \mathcal{S}$
such that $s=\lambda(|\xi|,s)$ to construct solutions to the
original problem (1.48)--(1.52).

\begin{lemma}There exists a unique $s\in \mathcal{S}$ so that
$\lambda(|\xi|,s)=\sqrt{-\mu( s)}>0$ and
\begin{equation}s=\lambda(|\xi|,s).\end{equation}
\end{lemma}
\hspace{-18pt}{\bf Proof.} By Lemma 3.5, there exists $s_\ast>0$
such that
\begin{equation}\mathcal{S}=\mu^{-1}((-\infty,0))=(0,s_\ast).\end{equation}
We define $\lambda=\sqrt{-\mu}$ on $\mathcal{S}$ and define the
function $\Phi:\ (0,s_\ast)\rightarrow(0,\infty)$ by
\begin{equation}\Phi(s)=s/\lambda(|\xi|,s),\end{equation}
which is continuous and strictly increasing in $s$. Moreover,
$\lim_{s\rightarrow0}\Phi(s)=0$ and $\lim_{s\rightarrow
s_\ast}\Phi(s)=+\infty$. Hence there is unique $s\in (0,s_\ast)$ so
that $\Phi(s)=1$, which gives (3.44). The proof of Lemma 3.7 is
completed.\hfill$\Box$\smallskip\smallskip

In view of Proposition 3.6, Lemma 3.7 and Lemma 3.4, we conclude Theorem 2.1.\hfill$\Box$

\subsection{Proof of Theorem 2.2}
The proof of Theorem 2.2 is similar to that of Theorem 2.1. The strategy is the same and the only difference is the energies
defined when using the variational method to construct solutions. We define the two energies related to the problem (1.56)--(1.60) by
\begin{eqnarray}&& E(\psi)
=\frac{1}{2}\int_{-1}^1s\mu(4|\xi|^2|\psi'|^2+||\xi|^2\psi+\psi''|^2)+|B|^2(|\xi|^4|\psi|^2+|\xi|^2|\psi'|^2)\,dx_2\nonumber
\\&&\qquad\quad\ -\frac{1}{2}|\xi|^2g[\rho]|\psi(0)|^2,
\\&& J(\psi)
=\frac{1}{2}\int_{-1}^1 \rho(|\xi|^2|\psi|^2+|\psi'|^2)\,dx_2,\end{eqnarray}
  which are both
well-defined on the space $H_0^2((-1,1))$.

We define
\begin{eqnarray}& E_0'(\psi)
&=\frac{1}{2}\int_{-1}^1|B|^2(|\xi|^2|\psi|^2+|\psi'|^2)\,dx_2-\frac{1}{2}g[\rho]|\psi(0)|^2.
\end{eqnarray}
As in the proof of Theorem 2.1, to prove Theorem 2.2 the first thing is to clarify the sign of $E_0'(\psi)$.  We have the following lemma.
\begin{lemma} We have the following three assertions:

$(i)$ If $|B|\ge |B|_c$, then for any $\psi$ we have
$E_0(\psi)\ge0$. Moreover,
\begin{equation}E_0'(\psi)\ge\frac{1}{2}\int_{-1}^1
(|B|^2-|B|_c^2)
|\psi'|^2+{|B|^2|\xi|^2|\psi|^2}\,dx_2.\end{equation}

$(ii)$ Fixed $|B|< |B|_c$ and $|\xi|\ge |\xi|_{hc}^B$, then for any
$\psi$  we have $E_0(\psi)\ge 0$.

$(iii)$ Fixed $|B|< |B|_c$ and $0<|\xi|<|\xi|_{hc}^B$, then there exists
$\psi$ such that $E_0'(\psi)<0$.

\end{lemma}
 \hspace{-18pt}{\bf Proof.} The proof follows by the definitions (2.4) and (2.6), which is similar to the that of Lemma 3.3.\hfill$\Box$\smallskip\smallskip

Once Lemma 3.8 is established, Theorem 2.2 follows similarly by the proof of Theorem 2.1.\hfill$\Box$

\subsection{Proof of Theorem 2.3}

In this subsection we assume  that $|B|<|B|_c$ to prove Theorem 2.3. Notice that in Lemma 3.7 the fixed point $s\in\mathcal{S}$ is unique, we may
write uniquely $\lambda(|\xi|)$ within $(|\xi|_{vc}^B,\infty)$ (resp. $(0,|\xi|_{hc}^B)$) for $\bar{B}$ is vertical (resp. $\bar{B}$ is horizontal),
while the
corresponding solution constructed in Theorem 2.1 $(iii)$ (resp. Theorem 2.2 $(iii)$) is written by $\psi_{|\xi|}$.

The continuity assertion follows in the same way as
in \cite[Proposition 3.9]{GT2}, it suffices to prove (2.7) and (2.8). To prove (2.7), we first derive the limit behavior. For this we take any
$|\xi|_n\in (|\xi|_{vc}^B,\infty)$  so that $|\xi|_n \rightarrow|\xi|_{vc}^B$,
then by Theorem 2.1 $(iii)$ there exist  functions
$\psi_{|\xi|_n}\in\mathcal{A}$ so that
\begin{equation}-\lambda^2(|\xi|_n)=E(\psi_{|\xi|_n})<0.\end{equation} Recalling the
expressions of energies (3.17), (3.18) and (3.19),  we have
\begin{equation}0<
\lambda^2(|\xi|_n)\le\frac{|\xi_n|^2g[\rho]}{2}\psi_{|\xi|_n}^2(0)-
\frac{|B|^2}{2}\int_{-1}^1(|\xi_n|^2|\psi_{|\xi|_n}|^2+|\psi_{|\xi|_n}''|^2)\,dx_2.
\end{equation}
We deduce from (3.51) and the fact
$\psi_{|\xi|_n}\in \mathcal{A}$ that
$\psi_{|\xi|_n}$ is uniformly bounded in $H_0^2 $ as in the proof of
Lemma 3.2. So, up to the extraction of a subsequence we have that
\begin{equation} \psi_{|\xi|_n} \rightarrow\widetilde{\psi} \hbox{ weakly in }H_0^2((-1,1)) \hbox{ and strongly in }H_0^1((-1,1))\end{equation}
as well as $\psi_{|\xi|_n}(0)\rightarrow\widetilde{\psi} (0)$. So
taking superior limit in (3.52) as $n\rightarrow\infty$ along the
subsequence, we have
\begin{equation}0\le
\limsup_{n\rightarrow\infty}\lambda^2(|\xi|_n)\le\frac{(|\xi|_{vc}^B)^2g[\rho]}{2}\widetilde{\psi}
^2(0)- \frac{|B|^2}{2}\int_{-1}^1((|\xi|_{vc}^B)^2|\widetilde{\psi}'
|^2+|\widetilde{\psi} ''|^2)\,dx_2\le 0.
\end{equation}
The last inequality above comes from Lemma 3.3 $(ii)$.
Since (3.54) holds for any such extracted subsequence, we deduce that
$\lim_{n\rightarrow\infty}\lambda^2(|\xi|_n)=0$ for the original
sequence $|\xi|_n$ as well. This proves $\lim_{|\xi|\rightarrow|\xi|_{vc}^B}\lambda(|\xi|)=0$.

Now we turn to prove the boundedness of $\lambda$. Using the fact that $-\lambda(|\xi|)^2=E(\psi_{|\xi|})$ and the
expressions of energies (3.17), (3.18) and (3.19), we have
\begin{eqnarray}&&\lambda(|\xi|)^2\le  \frac{1}{2} |\xi|^2g[\rho]|\psi_{|\xi|}(0)|^2-\frac{1}{2}\int_{-1}^1|B|^2(|\xi|^2|\psi_{|\xi|}'|^2+{|\psi_{|\xi|}''|^2})\,dx_2
.\end{eqnarray}
Notice additionally that $\psi_{|\xi|}\in \mathcal{A}$, we have
\begin{eqnarray} J(\psi_{|\xi|})
=\frac{1}{2}\int_{-1}^1 \rho(|\xi|^2|\psi_{|\xi|}|^2+|\psi_{|\xi|}'|^2)\,dx_2=1.\end{eqnarray}
Hence from (3.55)--(3.56) we can bound
$\lambda$ by
\begin{eqnarray}&\lambda(|\xi|)^2&\le \frac{g[\rho]}{2|B|^2\sqrt{\rho_+}}\left(|B|^2\int_{0}^1|\xi|^2|\psi_{|\xi|}'|^2\,dx_2\right)^\frac{1}{2}\left(\int_{0}^1\rho_+|\xi|^2|\psi_{|\xi|}|^2\,dx_2\right)^\frac{1}{2}\nonumber
\\&&\quad-\frac{1}{2}\int_{-1}^1|B|^2(|\xi|^2|\psi_{|\xi|}'|^2+{|\psi_{|\xi|}''|^2})\,dx_2\nonumber
\\&&\le \frac{2g[\rho]}{|B|^2\sqrt{\rho_+}}\int_{0}^1\rho_+|\xi|^2|\psi_{|\xi|}|^2\,dx_2\le \frac{4g[\rho]}{|B|^2\sqrt{\rho_+}}
.\end{eqnarray}
Hence, (2.7) follows.

Finally, to prove (2.8), the limit $\lim_{|\xi|\rightarrow|\xi|_{hc}^B}\lambda(|\xi|)=0$ can be proven in the same way as $\lim_{|\xi|\rightarrow|\xi|_{vc}^B}\lambda(|\xi|)=0$. It remains to prove $\lim_{|\xi|\rightarrow0}\lambda(|\xi|)=0$ and this follows from (3.55)--(3.56) again by
\begin{eqnarray}\lambda(|\xi|)^2\le \frac{|\xi|g[\rho]}{2\rho_+}\left(\int_{0}^1\rho_+|\xi|^2|\psi_{|\xi|}|^2\,dx_2\right)^\frac{1}{2}\left(\int_{0}^1\rho_+|\psi_{|\xi|}'|^2\,dx_2\right)^\frac{1}{2}\le \frac{|\xi|g[\rho]}{\rho_+}.\end{eqnarray}
The proof of Theorem 2.3 is completed.
\hfill$\Box$\smallskip\smallskip

\subsection{Proof of Theorem 2.4}

In this subsection we will do the energy estimates to prove Theorem 2.4. In this subsection and also in Theorem 2.4, $L^2,\ H^1$ denote the usual $L^p$ and Sobolev spaces on $\Omega$, and their norms are denoted by $\|\cdot\|_{L^2},\ \|\cdot\|_{H^1}$ respectively. We denote $C$ be generic constants depending only on the physical coefficients and $B$.

First, we prove $(i)$ and let $\bar{B}=(0,B)$ for the moment. To prove (2.10), differentiating the second equation in (1.35) with respect to time $t$ and eliminating the $\eta$ term by
using the first equation, we obtain
$$\left\{\begin{array}{ll}
\rho\partial_{tt}v +\nabla \partial_{t}q-\mu\Delta\partial_{t} v-|B|^2\partial_{22}^2v=0,
\\{\rm div} v={\rm div}\partial_{t}v=0,\end{array}\right.\eqno(3.59)$$
along with the jump and boundary conditions\setcounter{equation}{59}
\begin{eqnarray}&&\llbracket v\rrbracket=\llbracket \partial_{t}v\rrbracket=0,\quad\llbracket -\mu (D\partial_{t}v+D\partial_{t}v^T) + \partial_{t}qI\rrbracket e_2 =g[\rho]v_2e_2+|B|^2\llbracket\partial_{2}v\rrbracket,
\end{eqnarray}
\begin{equation}v_-(t,x_1,-1)=v_+(t,x_1,1)=\partial_{t}v_-(t,x_1,-1)=\partial_{t}v_+(t,x_1,1)=0.
\end{equation}

We regard the system (3.59)--(3.61) as one for $\partial_{t}v,\partial_{t}q,v$ and
we have the energy identity:
\begin{lemma}Let $\partial_{t}v,\partial_{t}q,v$ solve (3.59)--(3.61), then we have
\begin{eqnarray}
&&\frac{1}{2}\frac{d}{dt}\left(\int_\Omega\rho|\partial_{t}v|^2+|B|^2|\partial_2
v|^2\,dx-\int_{\mathbb{R}}g[\rho]|v_2(x_1,0)|^2\,dx_1\right)\nonumber
\\&&\quad+\int_\Omega\mu\left|
D\partial_{t}v+D\partial_{t}v^T\right|^2\,dx=0.
\end{eqnarray}
\end{lemma}
\hspace{-18pt}{\bf Proof.} Multiplying the equations $(3.59)_1$ by
$\partial_{t} v$ and integrate over $\Omega$, after integrating by
parts respectively in $\Omega_+$ and $\Omega_-$ and using the jump
and boundary conditions (3.60)--(3.61), by $(3.59)_2$, we obtain (3.62).\hfill$\Box$

\begin{lemma}Suppose $|B|> |B|_c$, then there exists constant $C^{-1}>0$ such that
\begin{eqnarray}\int_\Omega|B|^2|\partial_2
v|^2\,dx-\int_{\mathbb{R}}g[\rho]|v_2(x_1,0)|^2\,dx_1\ge C^{-1}(\|v\|_{L^2}^2+\|\partial_2v\|_{L^2}^2).\end{eqnarray}
\end{lemma}
\hspace{-18pt}{\bf Proof.} We represent the integral above as
\begin{eqnarray}&&\int_\Omega|B|^2|\partial_2
v|^2\,dx-\int_{\mathbb{R}}g[\rho]|v_2(x_1,0)|^2\,dx_1\nonumber
\\&&\quad=\int_{\mathbb{R}}\left(\int_{-1}^1|B|^2|\partial_2
v(x_1,x_2)|^2\,dx_2- g[\rho]|v_2(x_1,0)|^2\right)\,dx_1.\end{eqnarray}
By the definition (2.4) of the critical magnetic number $|B|_c$, we have that for any $x_1\in \mathbb{R}$,
\begin{eqnarray}&& \int_{-1}^1|B|^2|\partial_2
v(x_1,x_2)|^2\,dx_2- g[\rho]|v_2(x_1,0)|^2
\nonumber
\\&&\quad
\ge \int_{-1}^1\left(|B|^2|\partial_2
v_1(x_1,x_2)|^2+(|B|^2-|B|_c^2)|\partial_2
v_2(x_1,x_2)|^2\right)\,dx_2.\end{eqnarray}
Substituting (3.65) into (3.64), we have
\begin{eqnarray}\int_\Omega|B|^2|\partial_2
v|^2\,dx-\int_{\mathbb{R}}g[\rho]|v_2(x_1,0)|^2\,dx_1\ge C^{-1}\|\partial_2v\|_{L^2}^2.\end{eqnarray}
Hence, (3.63) follows from (3.66) and Poincar\'e inequality in the slab, since the boundary condition (3.61). The proof of Lemma 3.10 is completed. \hfill$\Box$\smallskip\smallskip

Now, we define
\begin{eqnarray}
&& \mathcal{E}(\partial_tv,v;t):=\int_\Omega\rho|\partial_tv|^2+|B|^2|\partial_2
v|^2\,dx-\int_{\mathbb{R}}g[\rho]|v_2(x_1,0)|^2\,dx_1,
\\&&\mathcal{D}(\partial_tv;t):=\int_\Omega\mu\left|
D\partial_tv+D\partial_tv^T\right|^2\,dx,\end{eqnarray}
then by Lemma 3.10, we have
\begin{eqnarray}
C^{-1}(\|\partial_tv\|_{L^2}^2+\|v\|_{L^2}^2+\|\partial_2v\|_{L^2}^2)\le \mathcal{E}(\partial_tv,v;t)\le C(\|\partial_tv\|_{L^2}^2+\|v\|_{L^2}^2+\|\partial_2v\|_{L^2}^2),\end{eqnarray}
and by Korn's inequality in the slab, since the boundary condition (3.61), we have
\begin{eqnarray}
C^{-1}\|\partial_tv\|_{H^1}^2\le \mathcal{D}(\partial_tv;t)\le C\|\partial_tv\|_{H^1}^2,\end{eqnarray}
Hence, integrating (3.62) directly in time, by (3.69)--(3.70), we obtain
\begin{eqnarray}&&
\|\partial_tv(t)\|_{L^2}^2+\|v(t)\|_{L^2}^2+\|\partial_2v(t)\|_{L^2}^2+\int_0^t\|\partial_tv(s)\|_{H^1}^2\,ds
\nonumber\\&&\quad\le C( \|\partial_tv(0)\|_{L^2}^2+\|v(0)\|_{L^2}^2+\|\partial_2v(0)\|_{L^2}^2).
\end{eqnarray}
This proves (2.10).

To prove (2.11), we notice that we have the boundary condition for $\eta$,
\begin{equation}\eta_-(t,x_1,-1)=\eta_+(t,x_1,1)=0,
\end{equation}
which is deduced from the first equation in (1.35) and the initial assumptions on $\eta^0$. Then $v,q,\eta$ satisfy the system (3.59)--(3.61) by replacing $\partial_{t}v,\partial_{t}q,v$ correspondingly.  Hence, the arguments for proving (2.10) also lead to (2.11).\smallskip\smallskip

Now we turn to prove $(ii)$ and hence let $\bar{B}=(B,0)$. The proof is similar to that of $(i)$. Similarly, to prove (2.12) we have the following system
$$\left\{\begin{array}{ll}
\rho\partial_{tt}v +\nabla\partial_t q-\mu\Delta\partial_t v-|B|^2\partial_{11}^2v=0,
\\{\rm div} v={\rm div}\partial_tv=0.\end{array}\right.\eqno(3.72)$$
along with the  jump and boundary conditions
\setcounter{equation}{72}
\begin{eqnarray}&&\llbracket v\rrbracket=\llbracket\partial_t v\rrbracket=0,\quad\llbracket -\mu (D\partial_tv+D\partial_tv^T) + \partial_tqI\rrbracket e_2 =g[\rho]v_2e_2,
\end{eqnarray}
\begin{equation} v_-(t,x_1,-1)= v_+(t,x_1,1)=\partial_tv_-(t,x_1,-1)=\partial_tv_+(t,x_1,1)=0.
\end{equation}

We have the  energy identity for the system (3.72)--(3.74):
\begin{lemma}Let $v,q$ solve (3.72)--(3.74), then we have
\begin{eqnarray}
&&\frac{1}{2}\frac{d}{dt}\left(\int_\Omega\rho|\partial_tv|^2+|B|^2|\partial_1
v|^2\,dx-\int_{\mathbb{R}}g[\rho]|v_2(x_1,0)|^2\,dx_1\right)\nonumber
\\&&\quad+\int_\Omega\mu\left|
D\partial_tv+D\partial_tv^T\right|^2\,dx=0.
\end{eqnarray}
\end{lemma}
\hspace{-18pt}{\bf Proof.} Multiplying the equations $(3.72)_1$ by
$v$ and integrate over $\Omega$, after integrating by
parts respectively in $\Omega_+$ and $\Omega_-$ and using the jump
and boundary conditions (3.73)--(3.74), by $(3.72)_2$, we obtain (3.75).\hfill$\Box$

\begin{lemma}Suppose $|B|> |B|_c$, then there exists constant $C^{-1}>0$ such that
\begin{eqnarray}\int_\Omega|B|^2|\partial_1
v|^2\,dx-\int_{\mathbb{R}}g[\rho]|v_2(x_1,0)|^2\,dx_1\ge C^{-1}(\|\partial_1v_1\|_{L^2}^2+\|v_2\|_{H^1}^2).\end{eqnarray}
\end{lemma}
\hspace{-18pt}{\bf Proof.} By $\partial_1v_1=-\partial_2v_2$, we can represent the integral in (3.76) as
\begin{eqnarray}&&\int_\Omega|B|^2|\partial_1
v|^2\,dx-\int_{\mathbb{R}}g[\rho]|v_2(x_1,0)|^2\,dx_1\nonumber
\\&&\quad=\int_\Omega|B|^2(|\partial_2
v_2|^2+|\partial_1
v_2|^2)\,dx-\int_{\mathbb{R}}g[\rho]|v_2(x_1,0)|^2\,dx_1\nonumber
\\&&\quad=\int_{\mathbb{R}}\left(\int_{-1}^1|B|^2(|\partial_2
v_2|^2+|\partial_1
v_2|^2)\,dx_2- g[\rho]|v_2(x_1,0)|^2\right)\,dx_1.\end{eqnarray}
By the definition (2.4) of  $|B|_c$, we have
\begin{eqnarray}&& \int_{-1}^1|B|^2|\partial_2
v_2(x_1,x_2)|^2\,dx_2- g[\rho]|v_2(x_1,0)|^2
\nonumber
\\&&\quad
\ge \int_{-1}^1 (|B|^2-|B|_c^2)|\partial_2
v_2(x_1,x_2)|^2 \,dx_2.\end{eqnarray}
Substituting (3.78) into (3.77), we have
\begin{eqnarray}\int_\Omega|B|^2|\partial_1
v|^2\,dx-\int_{\mathbb{R}}g[\rho]|v_2(x_1,0)|^2\,dx_1\ge C^{-1}(\|\partial_1v_1\|_{L^2}^2+\|v_2\|_{H^1}^2).\end{eqnarray}
Hence, (3.76) follows from (3.79) and Poincar\'e inequality in the slab, since the boundary condition (3.74). The proof of Lemma 3.12 is completed. \hfill$\Box$\smallskip\smallskip

Now, we define
\begin{eqnarray}
&& \mathcal{E}'(\partial_tv,v;t):=\int_\Omega\rho|\partial_tv|^2+|B|^2|\partial_1
v|^2\,dx-\int_{\mathbb{R}}g[\rho]|v_2(x_1,0)|^2\,dx_1,\end{eqnarray}
then by Lemma 3.12, we have
\begin{eqnarray}
C^{-1}(\|\partial_tv\|_{L^2}^2+\|\partial_1v_1\|_{L^2}^2+\|v_2\|_{H^1}^2)\le \mathcal{E}'(\partial_tv,v;t)\le C(\|\partial_tv\|_{L^2}^2+\|\partial_1v_1\|_{L^2}^2+\|v_2\|_{H^1}^2).\end{eqnarray}
Hence, integrating (3.75) directly in time, by (3.81) and Korn's inequality, we obtain
\begin{eqnarray}&&
\|\partial_tv(t)\|_{L^2}^2+\|\partial_1v_1(t)\|_{L^2}^2+\|v_2(t)\|_{H^1}^2+\int_0^t\|\partial_tv(s)\|_{H^1}^2\,ds
\nonumber\\&&\quad\le C( \|\partial_tv(0)\|_{L^2}^2+\|\partial_1v_1(0)\|_{L^2}^2+\|v_2(0)\|_{H^1}^2).
\end{eqnarray}
This proves (2.12).

To prove (2.13), by the initial condition and the first and third equations in (1.35),
\begin{equation}{\rm div}\eta=0.
\end{equation}
Hence, the situation is the same to the proof of (2.12) and hence along the similar lines we obtain (2.13). The proof of Theorem 2.4 is completed.
\hfill$\Box$

\section*{Acknowledgement}

The author would like to thank Yan Guo and Ian Tice for many helpful
discussions.   The author would like to express his gratitude for
the hospitality of the Division of Applied Mathematics at Brown
University during his visit, where this work was initiated.

\end{document}